\documentclass[11pt]{amsart}

\usepackage[margin=1in]{geometry}
\usepackage{amsmath,amssymb,amsthm,mathtools}
\usepackage{enumitem}
\usepackage{microtype}
\usepackage[hidelinks]{hyperref}
\usepackage{tikz}

\newtheorem{theorem}{Theorem}[section]
\newtheorem{proposition}[theorem]{Proposition}

\newtheorem{corollary}[theorem]{Corollary}
\newtheorem{remark}[theorem]{Remark}
\newtheorem{definition}[theorem]{Definition}
\newtheorem{question}[theorem]{Question}

\newcommand{\Mbar}{\overline{\mathcal M}}
\newcommand{\Mbarzn}{\overline{\mathcal M}_{0,n}}
\newcommand{\Mbargn}{\overline{\mathcal M}_{g,n}}

\newcommand{\RR}{\mathbb R}

\newcommand{\PP}{\mathbb{P}} 
\newcommand{\cM}{\mathcal{M}} 
\newcommand{\bcM}{\overline{\cM}}

\newcommand{\Pt}{\widetilde{P}}
\newcommand{\Ft}{\widetilde{F}}
\def\cM{\mathcal{M}}
\def\ZZ{\mathbb{Z}}
\def\CC{\mathbb{C}}
\def\LL{\mathbb{L}}
\def\lra{\longrightarrow}
\def\beq{\begin{equation}}
\def\eeq{\end{equation}}

\def\HH{\mathbb{H}}
\def\FF{\mathbb{F}}
\def\KV{\mathrm{K}(\mathbf{V}) }
\def\KVG{\mathrm{K}(\mathbf{V},G)}

\usepackage{pgfplots}
\pgfplotsset{compat=1.18}

\title{Stratified motivic invariants and bivariate deformations of Poincar\'e polynomials}
\author{Gergely B\'erczi and Young-Hoon Kiem}
\date{}
\address{Department of Mathematics, Aarhus University, Ny Munkegade 118, 8000 Aarhus, Denmark}
\email{gergely.berczi@math.au.dk}
\address{School of Mathematics, Korea Institute for Advanced Study, 85 Hoegiro, Dongdaemun-gu, Seoul 02455, Korea}
\email{kiem@kias.re.kr}

\begin{document}
\begin{abstract}
    For a stratified variety $X$ and a motivic invariant $\HH$, we consider the stratified invariant which interpolates the invariant of $X$ and that of its interior $X^\circ$. Based on observations in moduli theory, we introduce the notion of an echelon tower of stratified varieties and then prove explicit inductive formulae for the stratified invariants of the moduli spaces $\Mbarzn$ of stable curves of genus $0$ and the Fulton-MacPherson varieties $Y[n]$ for any smooth projective variety $Y$. Using these, we show that the mysterious bivariate deformations of the even degree Poincar\'e polynomials of $\bcM_{0,n+1}$ and $\PP^1[n]$ in \cite{BercziKiem2026} are nothing but stratified virtual Poincar\'e polynomials.    
\end{abstract}

\maketitle

\section{Introduction}

Let $\HH$ denote a motivic invariant such as the Hodge-Deligne polynomial or the virtual Poincar\'e polynomial. Let $X$ be a quasi-projective variety of dimension $m$ over $\CC$ with a stratification 
$$X=X_0\supset X_1\supset \cdots \supset X_m\supset X_{m+1}=\emptyset$$
by closed subvarieties such that $X_j^\circ=X_j-X_{j+1}$ is smooth of codimension $j$ or empty for each $j$. Then we can consider the \emph{stratified motivic invariant}
\beq\label{i3} \FF_X(y)=\sum_{j=0}^m y^j\HH_{X_j^\circ}\eeq 
which refines the motivic invariant of $X$ since $\FF_X(1)=\HH_X$. 
The first question that we want to address in this paper is the following.

\begin{question}
    How can we compute the stratified motivic invariants?
\end{question}

Important spaces in moduli theory often admit natural stratifications and come in families which are closely connected. For example, the moduli spaces of stable curves $\Mbargn$ are stratified by the dual graphs and are connected by the forgetful morphisms $\pi_n:\bcM_{g,n+1}\to \Mbargn$. Likewise, the Fulton-MacPherson compactifications $Y[n]$ of the configuration spaces of $n$ ordered distinct points in a smooth projective variety $Y$ are stratified by the dual graphs and connected by the forgetful morphisms $\pi_n:Y[n+1]\to Y[n]$.
As the forgetful morphisms arise from contracting a rational component if necessary after deleting the last marked point, the number of irreducible components reduces at most by 1. 
We can codify these observations as follows.

\begin{definition} [Definition \ref{m19}]
    A sequence of stratified spaces $X[n]$ for $n\ge n_0$ together with morphisms $\pi_n:X[n+1]\to X[n]$ are said to form an \emph{echelon tower} if $\pi_n(X[n+1]_j^\circ)\subset X[n]_j^\circ\cup X[n]_{j-1}^\circ.$ We call the echelon tower \emph{solid} if the restrictions 
    $$\pi_{n,j,0}:X[n+1]_j^\circ\cap \pi_n^{-1}(X[n]_j^\circ)\lra X[n]_j^\circ, \ \ \ 
\pi_{n,j,1}:X[n+1]_j^\circ\cap \pi_n^{-1}(X[n]_{j-1}^\circ)\lra X[n]_{j-1}^\circ $$
    of $\pi_n$ are motivically locally trivial. 
\end{definition}
For a solid echelon tower, we can compute the stratified invariant $\FF_{X[n]}$ by an inductive formula. 

\begin{theorem} [Proposition \ref{m15} and Theorem \ref{m14}] \label{i1}
    The moduli spaces $\Mbarzn$ of stable curves of genus $0$ form a solid echelon tower and the stratified invariant $\FF_{\Mbarzn}(y)$ is uniquely determined by the inductive formula
    \beq\label{i5} \FF_{\bcM_{0,n+1}} = (ny-n+1+\LL)\FF_{\Mbarzn}+y(y+\LL-1)\,\partial_y\FF_{\Mbarzn}, \quad \FF_{\bcM_{0,3}}=1\eeq
where $\LL=\HH_\CC$.
\end{theorem}
This is an effective way of computing the stratified invariants because the formula does not involve a sum over graphs or a sum over the product of $\bcM_{0,k}$ for all $k< n$ or inverting a complicated power series.  
We can easily derive previously known formulae for the Poincar\'e polynomial of $\Mbarzn$ from \eqref{i5}. 

\begin{theorem}[Theorem \ref{f8}] \label{i2}
    The Fulton-MacPherson varieties $Y[n]$ for a smooth projective variety $Y$ form a solid echelon tower and $\FF_{Y[n]}$ is uniquely determined by 
    \[\FF_{Y[n+1]}=\left(\HH_Y-n+ny\HH_{\PP^{m-1}}\right)\FF_{Y[n]}+ y\left(y\HH_{\PP^{m-1}}+ ({\LL^m}-1)\right)\partial_y\FF_{Y[n]},\quad \FF_{Y[1]}=\HH_Y.\]
\end{theorem}
From this, by specializing to $y=1$, we can derive previously known formulae for the Poincar\'e polynomials of $Y[n]$ including Manin's formulae \cite[(4.22) and (4.24)]{Manin1999} by elementary calculus.

\medskip

In \cite{BercziKiem2026}, it was proved that the even degree Poincar\'e polynomial of $\Mbarzn$ defined by 
\beq\label{i7}
P^{even}_{n}(t)=\sum_it^i\dim H^{2i}(\bcM_{0,n})\quad \text{or}\quad P^{even}_{n}(t^2)=P_{\Mbarzn}(t)\eeq 
has only real roots and hence the even degree Betti numbers form an ultra-log-concave sequence. In fact, the proof was assisted by AI which introduced, seemingly out of nowhere, a bivariate polynomial $F_n(y,t)$ such that $F_n(1,t)$ coincides with $P^{even}_{n+1}(t)$.
Obviously, the real curve in $\RR^2$ defined by $F_n(y,t)=0$ meets the line $y=1$ at the real roots of $P^{even}_{n+1}(t)$. By induction, it was shown that for a fixed $t=0^-$, $F_n(y,t)$ has $n-2$ distinct roots in the interval $(0,1)$ while for a sufficiently negative fixed $t$, $F_n(y,t)$ has $n-2$ distinct roots in the interval $(1,\infty)$. Since the curve is smooth by the implicit function theorem, it follows from a Sturm-Liouville argument and a degree count that $P^{even}_{n+1}(t)$ has $n-2$ simple negative roots. The same was proved for the Fulton-MacPherson variety $\PP^1[n]$. 

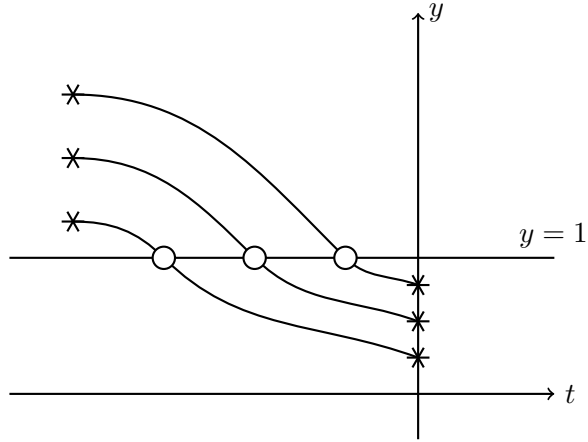
\begin{center}
\begin{tikzpicture}[
    scale=1.2,
    dot/.style={circle, fill=black, inner sep=1.5pt},
    arrow/.style={->, thick, >=stealth}
]

\tikzset{
    asterisk/.pic={
        \draw[thick] (-0.15,0) -- (0.15,0);
        \draw[thick] (-0.075,-0.13) -- (0.075,0.13);
        \draw[thick] (-0.075,0.13) -- (0.075,-0.13);
    }
}

\draw[->, thick] (-4.5, 0) -- (1.5, 0) node[right] {$t$};
\draw[->, thick] (0, -0.5) -- (0, 4.2) node[right] {$y$};

\draw[thick] (-4.5, 1.5) -- (1.5, 1.5);
\node[above right] at (1, 1.5) {$y=1$};

\draw[thick] (-3.8, 3.3) to[out=0, in=135] (-0.8, 1.5) to[out=-45, in=160] (0, 1.2);
\draw[thick] (-3.8, 2.6) to[out=0, in=135] (-1.8, 1.5) to[out=-45, in=160] (0, 0.8);
\draw[thick] (-3.8, 1.9) to[out=0, in=135] (-2.8, 1.5) to[out=-45, in=160] (0, 0.4);

\pic at (-3.8, 3.3) {asterisk};
\pic at (-3.8, 2.6) {asterisk};
\pic at (-3.8, 1.9) {asterisk};

\pic at (0, 1.2) {asterisk};
\pic at (0, 0.8) {asterisk};
\pic at (0, 0.4) {asterisk};

\filldraw[fill=white, thick] (-0.8, 1.5) circle (3.5pt);
\filldraw[fill=white, thick] (-1.8, 1.5) circle (3.5pt);
\filldraw[fill=white, thick] (-2.8, 1.5) circle (3.5pt);

\end{tikzpicture}
\end{center}
\medskip

Although the proof of the real-rootedness of the even degree Poincar\'e polynomial $P^{even}_{n}(t)$ of $\Mbarzn$ is complete, the geometric meaning of the key ingredient $F_n(y,t)$ has remained a mystery. So it is reasonable to ask the following. 
\begin{question}\label{i4}
    Is there a natural geometric meaning of the bivariate deformation $F_n(y,t)$ (resp. $\widetilde{F}_n(y,t)$) of the even degree Poincar\'e polynomial of $\bcM_{0,n+1}$ (resp. $\PP^1[n]$) in \cite{BercziKiem2026}?
\end{question} 

As an application of Theorems \ref{i1} and \ref{i2}, we obtain the following answer to Question \ref{i4}.
\begin{corollary} [Corollary \ref{m17} and Corollary \ref{f21}]\label{i6}
    (1) The bivariate deformation $F_n(y,t^2)$ of the Poincar\'e polynomial of $\bcM_{0,n+1}$ coincides with $y\FF_{\bcM_{0,n+1}}(y,t)$ where 
    $\HH$ is the virtual Poincar\'e polynomial. 

    (2) The bivariate deformation $\widetilde{F}_n(y,t^2)$ of the Poincar\'e polynomial of $\PP^1[n]$ in \cite{BercziKiem2026} is precisely the stratified virtual Poincar\'e polynomial $\FF_{\PP^1[n]}(y,t)$.  
\end{corollary} 
The $t^2$ appears in Corollary \ref{i6} from 
$P^{even}_{n}(t^2)=P_{\Mbarzn}(t)$ and the same for $\PP^1[n]$.

\medskip

It will be interesting to see if the ideas in this paper and \cite{BercziKiem2026} extend to Chow rings of matroids and give us the real-rootedness of the Poincar\'e polynomials of loopless matroids.

\bigskip

\section{Motivic invariants, stratified varieties and echelon towers}
In this section, we introduce the notions of motivic invariants, stratified varieties and echelon towers. 

\subsection{Motivic invariants of stratified varieties}
Let $\KV$ denote the Grothendieck ring of quasi-projective varieties over $\CC$. A \emph{motivic invariant} is a homomorphism 
$$\HH:\KV\lra R,\quad [X]\mapsto \HH_X$$
to an integral domain $R$ with unity such that $\HH_{[\mathrm{pt}]}=1$ and $\LL:=\HH_\CC$ invertible. 
In particular, we have 
\begin{enumerate}
    \item $\HH_X=\HH_Y+\HH_{X-Y}$ for a closed subvariety $Y\subset X$ and 
    \item $\HH_{X\times Z}=\HH_X \cdot \HH_Z$ for varieties $X$ and $Z$. 
\end{enumerate} 
As a consequence, if $X\to Y$ is a Zariski locally trivial morphism with fiber $Z$, we have $\HH_X=\HH_Y\cdot \HH_Z$. 
In what follows, we will use suitable extensions of the ring $R$, such as its field of fractions, whenever necessary.


A typical example of a motivic invariant is 
the Hodge-Deligne polynomial of a quasi-projective variety $X$ defined by  
$$\HH_X=\sum_{k,p,q}(-1)^ku^{p}v^q\dim H^{p,q;k}_c(X) \in \ZZ[u,v],\ \ \text{where } H^{p,q;k}(X)=\mathrm{gr}^p_F\mathrm{gr}^W_{p+q}H^k_c(X).$$
Here $\mathrm{gr}^p_F$ and $\mathrm{gr}^W_{p+q}$ denote the graded part by the Hodge and weight filtrations respectively of the mixed Hodge structure on $H^k_c(X)$. Letting $u=v=-t$, we get the virtual Poincar\'e polynomial 
$$P_X(t)=\sum_{k,p,q}(-1)^{k+p+q}t^{p+q}\dim H^{p,q;k}_c(X)\ \in \ \ZZ[t]$$
which coincides with the ordinary Poincar\'e polynomial $\sum_kt^k\dim H^k(X)$ when $X$ is smooth projective. 

\medskip

A \emph{stratified variety} is a variety $X$ together with a filtration
\beq\label{m1}
X=X_0\supset X_1\supset X_2\supset \cdots \supset X_{m}\supset X_{m+1}=\emptyset
\eeq
by closed subvarieties such that each $X_j^\circ=X_j-X_{j+1}$ is a smooth quasi-projective variety of codimension $j$ where $m=\dim X$. Let $X^\circ=X_0^\circ=X-X_1$.

\begin{definition}
    Let $\HH$ be a motivic invariant of varieties and $X$ be a variety with a  stratification \eqref{m1}. Then the \emph{stratified motivic invariant} of $X$ is defined as
    \beq\label{m2} \FF_{X}(y) =\sum_{j\ge 0} y^j \HH_{X_j^\circ} =\HH_{X^\circ}+y\HH_{X_1^\circ}+\cdots+y^{m-1}\HH_{X_{m-1}^\circ}+y^m\HH_{X_m}\ \in \ R[y].\eeq
\end{definition}

By definition, $\FF_X$ interpolates $\HH_{X^\circ}$ and $\HH_{X}$ because 
$\FF_X(0)=\HH_{X^\circ}$ and $\FF_X(1)=\HH_X$.
Obviously $\FF_{X}$ has more information than $\HH_X$ but we will see below that for echelon towers, the computation of $\FF_{X}$ may be easier than that of $\HH_X$ by an inductive formula. 

\begin{remark} There is an equivariant version of $\FF_X$ as follows. 
For a group $G$, let $\KVG$ denote the Grothendieck group of $G$-varieties. A \emph{$G$-equivariant motivic invariant} is a homomorphism 
$$\HH^G:\KVG\lra R, \quad [X]\mapsto \HH^G_X$$
to a commutative ring $R$ with unity such that $\HH^G_{\mathrm{pt}}=1$. 
For a $G$-invariant stratification \eqref{m1} on a $G$-variety $X$, we can define the \emph{equivariant stratified motivic invariant} $\FF^G_{X}$ by \eqref{m2} with $\HH$ replaced by $\HH^G$. 
\end{remark}

\subsection{Echelon towers of stratified varieties}
Many interesting moduli spaces in algebraic geometry come in a sequence like the moduli space $\Mbarzn$ of stable curves of genus $0$, the Fulton-MacPherson space $Y[n]$ or the stable map space $\Mbargn(Y,d)$ for a smooth projective variety $Y$.   
These spaces all admit natural stratifications and satisfy a nice property that we codify as follows. 

\begin{definition}\label{m19} An \emph{echelon tower} of stratified varieties consists of 
    a sequence of stratified varieties
        $$X[n]=X[n]_0\supset X[n]_1\supset X[n]_2\supset \cdots\ \ \text{ for each  }n\ge n_0$$ and morphisms $\pi_n:X[n+1]\to X[n]$ for $n\ge n_0$
    such that 
        $$\pi_{n}(X[n+1]_j^\circ)\subset X[n]_j^\circ\cup X[n]_{j-1}^\circ.$$
An echelon tower is called \emph{solid} if the morphisms   
\beq\label{m3}\pi_{n,j,0}:X[n+1]_j^\circ\cap \pi_n^{-1}(X[n]_j^\circ)\lra X[n]_j^\circ, \ \ 
\pi_{n,j,1}:X[n+1]_j^\circ\cap \pi_n^{-1}(X[n]_{j-1}^\circ)\lra X[n]_{j-1}^\circ\eeq 
induced by $\pi_n$ are motivically locally trivial and the fibers have constant motivic invariants $A_{n,j}$ and $B_{n,j}\in R$ respectively.
\end{definition}

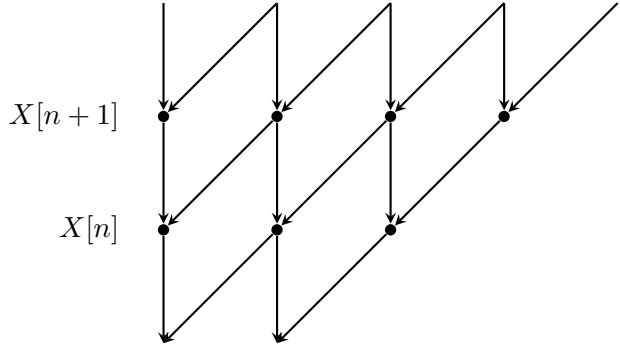
\begin{figure}
\centering
\begin{tikzpicture}[
    scale=1.5,
    dot/.style={circle, fill=black, inner sep=1.5pt},
    arrow/.style={->, thick, >=stealth}
]

\node[left] at (-0.3, 1) { $X[n+1]$};
\node[left] at (-0.3, 0) { $X[n]$};

\node[dot] (n1_1) at (0, 1) {};
\node[dot] (n1_2) at (1, 1) {};
\node[dot] (n1_3) at (2, 1) {};
\node[dot] (n1_4) at (3, 1) {};

\node[dot] (n0_1) at (0, 0) {};
\node[dot] (n0_2) at (1, 0) {};
\node[dot] (n0_3) at (2, 0) {};

\foreach \i in {1,2,3,4} {
    \draw[arrow] (n1_\i |- 0, 2) -- (n1_\i);
    \draw[arrow] ([xshift=1cm, yshift=1cm]n1_\i.center) -- (n1_\i);
}

\draw[arrow] (n1_1) -- (n0_1);
\draw[arrow] (n1_2) -- (n0_2);
\draw[arrow] (n1_3) -- (n0_3);

\draw[arrow] (n1_2) -- (n0_1);
\draw[arrow] (n1_3) -- (n0_2);
\draw[arrow] (n1_4) -- (n0_3);

\draw[arrow] (n0_1) -- (0, -1);
\draw[arrow] (n0_2) -- (1, -1);

\draw[arrow] (n0_3) -- (1, -1);
\draw[arrow] (n0_2) -- (0, -1);

\end{tikzpicture}
\caption{The dots represent $X[n]_j^\circ$. The vertical arrows are $\pi_{n,j,0}$ and the slanted arrows are $\pi_{n,j,1}$.}
\end{figure}
\medskip

Here we say that a morphism $f:U\to Y$ is \emph{motivially locally trivial} if it factors as
$$U\hookrightarrow X \stackrel{g}{\longrightarrow} Y$$
such that $U$ is open in $X$ and that the morphisms $g$ and $g|_{X-U}$ are Zariski locally trivial. 
Then for any motivic invariant $\HH$, we have $\HH_U=\HH_F\cdot\HH_Y$ where $F$ denotes the fiber of $f$ if $Y$ is connected. 

By \eqref{m3}, we have
\beq\label{m4}
\HH_{X[n+1]_j^\circ}=A_{n,j}\HH_{X[n]_j^\circ}+ B_{n,j}\HH_{X[n]_{j-1}^\circ}.\eeq
If $A_{n,j}$ and $B_{n,j}$ are polynomials in $j$, then we will see that \eqref{m4} induces a recursive formula for $\FF_{X[n]}$, which is an effective way of computing $\FF_{X[n]}$ and $\HH_{X[n]}=\FF_{X[n]}(1)$.

\medskip

In the subsequent sections, we will use \eqref{m4} to deduce formulae for stratified motivic invariants of 
the moduli spaces of stable curves of genus $0$ and the Fulton-MacPherson varieties.

\bigskip

\section{Moduli spaces of stable curves of genus 0}
A fundamental example of an echelon tower is the moduli spaces $\Mbarzn$ of stable curves of genus $0$ with $n$ marked points for $n\ge 3$. Recall that $\Mbarzn$ consists of connected trees of projective lines with $n$ ordered smooth distinct marked points, each of whose irreducible components contains at least three special (nodal or marked) points. It compactifies the moduli space $\cM_{0,n}$ of $n$ distinct ordered points in $\PP^1$ up to projective equivalence. 

To each $C\in \Mbarzn$, we can associate its dual graph whose vertices (resp. edges, ordered legs) correspond to irreducible components (resp. nodes, marked points). For a graph $\gamma$, let $e_\gamma$ denote the number of edges. Since the genus is 0, the dual graph is always a tree and the number of vertices is $e_\gamma+1$. 
For each tree $\gamma$ with $n$ ordered legs, let $\cM_{0,n,\gamma}$ be the locus of $C\in \Mbarzn$ whose dual graph is $\gamma$. Then we have
$$\cM_{0,n,\gamma}=\prod_{v\in \gamma}\cM_{0,k(v)}$$
where $k(v)$ denotes the sum of the number of nodes and the number of marked points in 
the irreducible component corresponding to a vertex $v$.  In particular, $\cM_{0,n,\gamma}$ is a smooth subvariety of codimension $e_\gamma$ in $\Mbarzn$.
Let 
\beq\label{m6}
\cM_{0,n,j}=\bigsqcup_{\gamma:e_\gamma=j}\cM_{0,n,\gamma}.\eeq
Then we have a stratification 
\beq\label{m9}
\Mbarzn=\bcM_{0,n,0}\supset \bcM_{0,n,1}\supset \bcM_{0,n,2}\supset\cdots\supset \bcM_{0,n,n-3}\supset \bcM_{0,n,n-2}=\emptyset\eeq
with $\cM_{0,n}=\cM_{0,n,0}$ and the stratified motivic invariant of $\Mbarzn$ is defined as
\beq\label{m8}
\FF_{\Mbarzn}(y) =\sum_{j=0}^{n-3} y^j\HH_{\cM_{0,n,j}}. 
\eeq

For $n\ge 3$, we have the morphism $$\pi_n:\bcM_{0,n+1}\lra \Mbarzn$$ which forgets the last marked point. 
For $C\in \bcM_{0,n+1}$, $\pi_n(C)$ is obtained by contracting the component in $C$ containing the last marked point if it has only three special points. Otherwise, the underlying curves of $C$ and $\pi_n(C)$ are the same. 
Thus we find that 
$$\cM_{0,n+1,j} =  \left(\cM_{0,n+1,j}\cap \pi_n^{-1}(\cM_{0,n,j})\right) \cup \left(\cM_{0,n+1,j}\cap \pi_n^{-1}(\cM_{0,n,j-1})\right)$$
and obtain morphisms
\beq\label{m5} \pi_{n,j,0}: \cM_{0,n+1,j}\cap \pi_n^{-1}(\cM_{0,n,j})\lra \cM_{0,n,j},\eeq
\[ \pi_{n,j,1}: \cM_{0,n+1,j}\cap \pi_n^{-1}(\cM_{0,n,j-1})\lra \cM_{0,n,j-1}.\]

By fixing the first three marked points at $0,1,\infty$, the open moduli space $\cM_{0,n}$ can be identified with 
$$(\CC-\{0,1\})^{n-3}-\bigcup_{i<j}\{(p_4,\cdots,p_n)\,|\,p_i=p_j\}$$
and $\pi_n:\cM_{0,n+1}\to \cM_{0,n}$ is the restriction of the projection $\cM_{0,n}\times (\CC-\{0,1\})\to \cM_{0,n}$ to the complement of $n-3$ sections.  
Hence, $\pi_n:\cM_{0,n+1}\to \cM_{0,n}$ is motivically locally trivial and 
$$\HH_{\cM_{0,n+1}}=\HH_{\cM_{0,n}} (\LL-n+1),\quad \LL=\HH_\CC.$$
By induction with $\cM_{0,3}=\mathrm{pt}$, we find that 
$$\HH_{\cM_{0,n+1}}=\prod_{j=2}^{n-1}(\LL -j)=\frac{n!}{\LL(\LL-1)}\binom{\LL}{n}.$$
Likewise, by \eqref{m6}, we can decompose \eqref{m5} into the unions of morphisms over $\cM_{0,n,\gamma}$ for $e_\gamma=j$ or $j-1$ and find that the morphisms in \eqref{m5} are motivically locally trivial. The fiber over $C\in \cM_{0,n,j}$ of $\pi_{n,j,0}$ is $C$ minus the $n+j$ special points and has the motivic invariant $$(j+1)\LL -j -n+1.$$ 
The fiber over $C\in \cM_{0,n,j-1}$ of $\pi_{n,j,1}$ consists of $n+j-1$ special points. 
So we proved the following.
\begin{proposition}\label{m15} The moduli spaces $\bcM_{0,n}$ for $n\ge 3$ form a solid echelon tower of stratified varieties and the motivic invariants of strata satisfy  
\beq\label{m7}
\HH_{\cM_{0,n+1,j}}=\left((j+1)\LL -j -n+1\right) \HH_{\cM_{0,n,j}} + (n+j-1)\HH_{\cM_{0,n,j-1}}.
\eeq
\end{proposition}

By \eqref{m7}, we have  
$$\FF_{\bcM_{0,n+1}}=\sum_j y^{j}\HH_{\cM_{0,n+1,j}}=\sum_j(j(\LL-1)+\LL -n+1)y^{j}\HH_{\cM_{0,n,j}}
+\sum_j (n+j-1)y^j\HH_{\cM_{0,n,j-1}}$$
$$=(\LL-1)y\sum_jjy^{j-1}\HH_{\cM_{0,n,j}}+(\LL -n+1)\sum_jy^j\HH_{\cM_{0,n,j}}$$ 
$$+ny\sum_jy^{j-1}\HH_{\cM_{0,n,j-1}}+y^2\sum_j (j-1)y^{j-2}\HH_{\cM_{0,n,j-1}}$$
\[
=(ny-n+1+\LL)\FF_{\Mbarzn}+y(y+\LL-1)\,\partial_y\FF_{\Mbarzn}.\]
So we obtain the following.
\begin{theorem}\label{m14}
    The stratified motivic invariants $\FF_{\Mbarzn}(y)$ of $\Mbarzn$ are uniquely determined by the inductive formula
    \beq\label{m11} \FF_{\bcM_{0,n+1}} = (ny-n+1+\LL)\FF_{\Mbarzn}+y(y+\LL-1)\,\partial_y\FF_{\Mbarzn}, \quad \FF_{\bcM_{0,3}}=1.\eeq
\end{theorem}
For example, $\FF_{\bcM_{0,4}}=(\LL-2)+3y$ and 
$$\FF_{\bcM_{0,5}}=(\LL-2)(\LL-3)+(10\LL-20)y+15y^2.$$

Note that the right hand side of \eqref{m11} equals $$((n-1)y-n+2)\FF_{\Mbarzn}+(y+\LL-1)\partial_y(y\FF_{\Mbarzn})$$
since $\partial_y(y \FF_{\Mbarzn})=\FF_{\Mbarzn}+y\partial_y\FF_{\Mbarzn}$. Multiplying $y$, we obtain 
$$y\FF_{\bcM_{0,n+1}}=((n-1)y-n+2)y\FF_{\Mbarzn}+y(y+\LL-1)\partial_y(y\FF_{\Mbarzn}).$$
Let $F_{n}(y)=y\FF_{\bcM_{0,n+1}}(y)$. Then we have 
\beq\label{m10}
F_{n+1}(y)=(ny-n+1)F_n(y)+y(y+\LL-1)\partial_yF_n(y), \quad F_2(y)=y\FF_{\bcM_{0,3}}=y. 
\eeq
Note that \eqref{m10} is precisely the defining inductive formula for the bivariate deformation $F_n(y,t)$ of the even degree Poincar\'e polynomial of $\bcM_{0,n+1}$ in \cite[Definition 2.2]{BercziKiem2026} when $\HH$ is the virtual Poincar\'e polynomial. So we proved the following.
\begin{corollary}\label{m17}
    The bivariate deformation $F_n(y,t)$ of the even degree Poincar\'e polynomial 
    of $\bcM_{0,n+1}$ in \cite{BercziKiem2026} is uniquely determined by   $F_n(y,t^2)=y\FF_{\bcM_{0,n+1}}(y)$ when $\HH$ is the virtual Poincar\'e polynomial.   
\end{corollary}
Consequently, the coefficient of $y^{j+1}$ in $F_n(y,t)$ is the virtual Poincar\'e polynomial of $\cM_{0,n+1,j}$:
\beq\label{m12}
[y^{j+1}]F_n(y,t^2)=P_{\cM_{0,n+1},j}({t}).\eeq
Thus Corollary \ref{m17} provides us with a geometric meaning of the coefficients of $F_n(y,t)$. 

See \cite{BercziKiem2026} 
for the generating function $$\Phi(x,y)=x+\sum_{n\ge 2}\frac{x^n}{n!}y\FF_{\bcM_{0,n+1}}(y),$$ the partial differential equation satisfied by $\Phi$, and a proof of the real-rootedness of the Poincar\'e polynomial of $\Mbarzn$ as an application.
We can recover all the previously known formulae for the Poincar\'e polynomial of $\Mbarzn$ from \eqref{m7} or \eqref{m11} by using the computation in \cite{BercziKiem2026}. 

\begin{remark}\label{m18}
    The moduli spaces $\Mbargn$ also form an echelon tower which is not solid. Everything holds except that the morphism $$\pi_{n,j,0}:\cM_{g,n+1,j}\cap \pi_n^{-1}(\cM_{g,n,j})\lra \cM_{g,n,j}$$
    is not motivically locally trivial where $\pi_n:\bcM_{g,n+1}\to \Mbargn$ denotes the forgetful morphism. Likewise, the moduli spaces of stable maps $\bcM_{g,n}(Y,d)$ for a convex smooth projective variety $Y$ form an echelon tower that is not solid. 
\end{remark}

\bigskip

\section{Stratified motivic invariants of Fulton-MacPherson varieties}

In this section, we fix a motivic invariant $\HH$ and determine the stratified motivic invariants of the Fulton-MacPherson varieties. 

For a smooth projective variety $Y$ of dimension $m$, let $Y[n]$ denote the Fulton-MacPherson compactification of the configuration space 
$$Y[n]^\circ=Y^n-\bigcup_{i<j}\Delta_{ij},\quad \Delta_{ij}=\{(p_1,\cdots,p_n)\,|\,p_i=p_j\}$$
of $n$ ordered distinct points in $Y$ (cf. \cite{FultonMacpherson1994}). 
Each element in $Y[n]$ parameterizes a variety $Y'$, obtained from $Y$ by repeatedly blowing up at a smooth point and gluing $\PP^m$ along the exceptional divisor, together with $n$ smooth distinct marked points $p_i$ such that the automorphism group is trivial. 

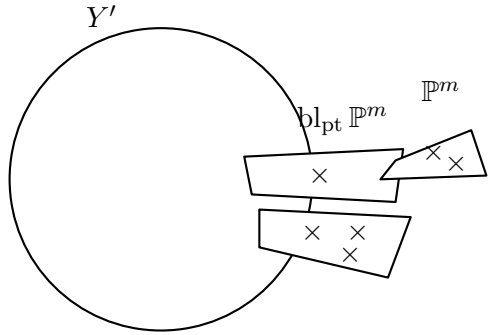
\begin{center}
\begin{tikzpicture}[line cap=round, line join=round, >=stealth]

  \draw[thick] (0,0) circle (2cm);
  \node at (110:2.3) {$Y'$};


  \filldraw[fill=white, thick] (1.3, -0.4) -- (3.3, -0.5) -- (3.0, -1.3) -- (1.3, -0.9) -- cycle;
  
  \node at (2.0, -0.7) {$\times$};
  \node at (2.6, -0.7) {$\times$};
  \node at (2.5, -1.0) {$\times$};

  \filldraw[fill=white, thick] (1.1, 0.3) -- (3.2, 0.4) -- (3.1, -0.3) -- (1.2, -0.2) -- cycle;
  \node[anchor=south] at (2.4, 0.45) {$\operatorname{bl}_{\operatorname{pt}} \mathbb{P}^m$};
  
  \node at (2.1, 0.05) {$\times$};

  \filldraw[fill=white, thick] (2.9, 0.0) -- (3.1, 0.25) -- (4.1, 0.65) -- (4.3, 0.05) -- cycle;
  \node[anchor=south] at (3.7, 0.9) {$\mathbb{P}^m$};

  \node at (3.6, 0.35) {$\times$};
  \node at (3.9, 0.2) {$\times$};

\end{tikzpicture}
\end{center}

For $(Y',p_1,\cdots, p_n)\in Y[n]$, let $\gamma_{Y'}$ be its dual graph whose vertices (resp. edges, legs) correspond to the irreducible components (resp. intersections $\PP^{m-1}$ of irreducible components, marked points) in $Y'$. There is a distinguished vertex representing (a blowup of) $Y$, called the distinguished component. The graph $\gamma_{Y'}$ is always a tree. 

For a tree $\gamma$ with $n$ legs and a distinguished vertex, let $Y[n]^\circ_\gamma$ be the subset of $Y[n]$ consisting of objects whose dual graph is $\gamma$.  
Then $Y[n]$ is stratified by  
$$Y[n]_j=\bigsqcup_{\gamma: e_\gamma\ge j} Y[n]^\circ_\gamma, \quad Y[n]^\circ_j=\bigsqcup_{\gamma: e_\gamma= j} Y[n]^\circ_\gamma$$
where $e_\gamma$ is the number of edges in $\gamma$.

Let $\pi_n:Y[n+1]\to Y[n]$ denote the morphism forgetting the last marked point $p_{n+1}$. If the irreducible component containing $p_{n+1}$ stays stable (i.e. the autormorphism group is trivial) after deleting it, the underlying variety remains the same. If not, we contract the irreducible component. Therefore, we have   
$$\pi_n\left( Y[n+1]^\circ_j\right) \subset Y[n]^\circ_j\sqcup Y[n]^\circ_{j-1}$$
so that 
\beq\label{f2}
Y[n+1]_j^\circ=\left(Y[n+1]_j^\circ\cap \pi_n^{-1}(Y[n]^\circ_j) \right)\sqcup \left(Y[n+1]_j^\circ\cap \pi_n^{-1}(Y[n]^\circ_{j-1}) \right).
\eeq 
Then $\pi_n$ induces morphisms 
\beq\label{f16}
\pi_{n,j,0}:Y[n+1]_j^\circ\cap \pi_n^{-1}(Y[n]^\circ_j) \lra Y[n]^\circ_j, 
\quad \pi_{n,j,1}:Y[n+1]_j^\circ\cap \pi_n^{-1}(Y[n]^\circ_{j-1}) \lra Y[n]^\circ_{j-1}\eeq 
As in the previous section, we can further decompose these into morphisms over $Y[n]^\circ_\gamma$ where $\gamma$ runs over trees with $e_\gamma=j, j-1$.  
The fiber of $\pi_{n,j,0}$ over $(Y',p_1,\cdots,p_n)$ is the complement in $Y'$ of  the marked points and the intersections of irreducible components in $Y'$. It is the same as the disjoint union of $Y$ and $j$ copies of $\PP^{m}$ minus $j$ points (blowup centers), $j$ copies of $\PP^{m-1}$ (intersections of irreducible components) and $n$ marked points. 
We thus find that $\pi_{n,j,0}$ is motivically locally trivial whose fibers have the motivic invariant
\beq\label{f4}
\HH_Y-j+j\HH_{\CC^m}-n=\HH_Y-n+j(\LL^m-1).
\eeq
The fiber of $\pi_{n,j,1}$ over $(Y',p_1,\cdots,p_n)$ consists of $(\widetilde{Y}',p_1,\cdots,p_{n+1} )\in Y[n+1]_j^\circ$ 
where $\widetilde{Y}'\to Y'$ contracts the irreducible component containing $p_{n+1}$ which has exactly one more marked point or one irreducible component connected in the direction away from the distinguished component.  
The contracted component is mapped to a marked point by the contraction $\widetilde{Y}'\to Y'$ in the former case and to the intersecting $\PP^{m-1}$ in the latter case. The automorphism group of $\PP^{m}$ fixing $\PP^{m-1}\cup\mathrm{pt}$ pointwisely or that of $\mathrm{bl}_{\mathrm{pt}}\PP^m$ fixing $\PP^{m-1}\cup E$ pointwisely is $\CC^*$ where $E$ denotes the exceptional divisor and so we find that 
$\pi_{n,j,1}^{-1}(Y',p_1,\cdots,p_n)$ is the union of 
of $n+j-1$ copies of $\PP^{m-1}=(\PP^m-\PP^{m-1}\cup \mathrm{pt})/\CC^*=(\mathrm{bl}_{\mathrm{pt}}\PP^m-\PP^{m-1}\cup E)/\CC^*$, each corresponding to a leg or an edge of the dual graph of $(Y',p_1,\cdots,p_n)$. 
We thus find that $\pi_{n,j,1}$ is also motivically locally trivial and 
the fibers have the motivic invariant 
\beq\label{f4'}
(n+j-1)\,\HH_{\PP^{m-1}}.
\eeq
Therefore we have
\beq\label{f5}
\HH_{Y[n+1]_j^\circ}= \left( \HH_Y-n+j(\LL^m-1) \right) \HH_{Y[n]^\circ_j} +
(n+j-1) \HH_{\PP^{m-1}} \HH_{Y[n]^\circ_{j-1}} 
\eeq
which immediately implies 
the following inductive formula for the Fulton-MacPherson space $Y[n]$ by a computation similar to that for \eqref{m10}.
\begin{theorem}\label{f8} The Fulton-MacPherson varieties $Y[n]$ form a solid echelon tower of stratified varieties.  
The stratified motivic invariant 
$\FF_{Y[n]}=\sum_j y^j\HH_{Y[n]_j^\circ}$ of $Y[n]$  is uniquely determined by the recursive formula
\beq\label{f6} 
\FF_{Y[n+1]}=\left(\HH_Y-n+ny\HH_{\PP^{m-1}}\right)\FF_{Y[n]}+ y\left(y\HH_{\PP^{m-1}}+ ({\LL^m}-1)\right)\partial_y\FF_{Y[n]}\eeq
with $\FF_{Y[1]}=\HH_Y$. 
\end{theorem}
Since $\FF_{Y[n]}(1)$ is the Hodge polynomial of ${Y[n]}$ when $\HH$ is the Hodge-Deligne polynomial, \eqref{f6} gives us an effective way of computing the Hodge and Poincar\'e polynomials of the Fulton-MacPherson variety $Y[n]$.

For example, when $Y=\PP^1$, \eqref{f6} is 
\beq\label{f17}
\FF_{\PP^1[n+1]}=(\LL+1-n+ny)\FF_{\PP^1[n]} + y(y+\LL-1)\partial_y\FF_{\PP^1[n]}.\eeq

\medskip

\noindent\textbf{Example.}
Since $Y[1]=Y$ is smooth, $\FF_{Y}=\HH_Y$. 
For $n=1$, \eqref{f6} is 
$$\FF_{Y[2]}=(\HH_Y-1+y\HH_{\PP^{m-1}})\HH_Y$$ which also follows from the blowup construction $Y[2]=\mathrm{bl}_Y(Y\times Y)$. 

For $n=2$, \eqref{f6} gives
$$\FF_{Y[3]}=(\HH_Y-2+2y\HH_{\PP^{m-1}})(\HH_Y-1+y\HH_{\PP^{m-1}})\HH_Y + \left(y^2\HH_{\PP^{m-1}}+ y(\LL^m-1)\right)\HH_{\PP^{m-1}}\HH_Y$$
$$=\HH_Y(\HH_Y-1)(\HH_Y-2)
+y(\LL^m+3\HH_Y-5)\HH_{\PP^{d-1}}\HH_Y
+3y^2\HH_{\PP^{d-1}}^2\HH_Y.$$
Letting $y=1$, we find that the Hodge polynomial of $Y[3]$ equals
$$\HH_Y^3+\HH_Y(\HH_{\PP^{2m-1}}-1)+3(\HH_Y^2+\HH_Y(\HH_{\PP^{m-1}}-1))(\HH_{\PP^{m-1}}-1)$$
which can be also read off from the blowup construction of $Y[3]$. In fact, $Y[3]\to Y^3$ is the composition of two blowups, first along the small diagonal $\Delta_{123}\cong Y$ and then along the union of the proper transforms of $\Delta_{12}$, $\Delta_{23}$ and $\Delta_{13}\cong Y^2$. 

\medskip

It is useful to define the generating function of all $\FF_{Y[n]}$. 
\begin{definition}
\beq\label{f18} \Psi(x,y)=1+\sum_{n\ge 1}\frac{x^n}{n!}\FF_{Y[n]}(y)\  \in \ R[y][\![x]\!].\eeq
\end{definition}
Then it is easy to check that \eqref{f6} is equivalent to the following partial differential equation.
\begin{proposition}
\beq\label{f7}
\HH_Y\Psi=(1+x(1-y\HH_{\PP^{m-1}}))\partial_x\Psi - y(y\HH_{\PP^{m-1}}+{\LL^m}-1)\partial_y\Psi.
\eeq
Moreover, $\Psi$ is uniquely determined by \eqref{f7} and $\Psi|_{x=0}=1$. 
\end{proposition}

When $Y=\PP^1$, \eqref{f7} gives 
$$(\LL+1)\Psi=(1-(y-1)x)\partial_x\Psi-y(y+\LL-1)\partial_y\Psi$$
and hence we find that \eqref{f18} coincides with the generating function $\Psi$ in \cite[\S5]{BercziKiem2026} with $t$ replaced by $\LL$. 
Since the coefficient of $\frac{x^n}{n!}$ in $\Psi$ is the bivariate deformation of the Poincar\'e polynomial of $\PP^1[n]$, we obtain the following.
\begin{corollary}\label{f21}
  The bivariate deformation $\widetilde{F}_n(y,t^2)$ of the Poincar\'e polynomial of $\PP^1[n]$ in \cite{BercziKiem2026} is precisely the stratified virtual Poincar\'e polynomial $\FF_{\PP^1[n]}(y,t)$.   
\end{corollary}

We end this paper with a new proof of Manin's formulae for the Poincar\'e polynomial of the Fulton-MacPherson variety $Y[n]$ in \cite{Manin1999}, from \eqref{f7}.

It is useful to define $\Phi\in R[y][\![x]\!]$ as the unique series satisfying
$$\Psi=(1+\Phi)^{\HH_Y}$$
when $\HH_Y\ne 0$.  
Then \eqref{f7} is equivalent to 
\beq\label{f9}
1+\Phi = (1+x(1-y\HH_{\PP^{m-1}}))\partial_x\Phi - y(y\HH_{\PP^{m-1}}+{\LL^m}-1)\partial_y\Phi.
\eeq

As in \cite{BercziKiem2026}, we can solve this partial differential equation along $y=1$ by using the characteristic equations: 
\beq\label{f11} \frac{dx}{ds}=1+x(1-y\HH_{\PP^{m-1}}),\quad \frac{dy}{ds}=- y(y\HH_{\PP^{m-1}}+{\LL^m}-1),\quad \frac{d\Phi}{ds}=1+\Phi.\eeq
Let $x=0$ at $s=0$ so that $\Phi=0$ at $s=0$ as well. Then the last equation gives us 
\beq\label{f10}
\Phi=e^s-1=u-1.\eeq
We can also solve the second equation in \eqref{f11} as
\beq\label{f12}
\frac{y}{y+\LL-1} = \frac{1}{\LL} \left( \frac{z}{u} \right)^{\LL^m-1}\eeq
where $z$ is the value of $u=e^s$ when $y(s)=1$, i.e. $y(\log z)=1$. 

From the first equation in \eqref{f11}, we have
\beq\label{f13}
\frac{dx}{du} + \frac{y\HH_{\PP^{m-1}} -1 }{u} x =\frac{1}{u}\eeq 
From \eqref{f12}, we find that 
$$y\HH_{\PP^{m-1}} -1  = \frac{\LL^m\left( \frac{z}{u} \right)^{\LL^m-1}-\LL}{\LL-\left( \frac{z}{u}\right) ^{\LL^m-1}}.$$
We thus have 
\beq\label{f14} \exp \int \frac{y\HH_{\PP^{m-1}} -1}{u}  \,du = \frac{z}{u}\left( \LL - \left(\frac{z}{u}\right)^{\LL^m-1}\right).\eeq
By multiplying \eqref{f14} to \eqref{f13} and integrating from $u=1$ to $u=z$ while remembering $x=0$ when $u=1$, we obtain
$$(\LL-1)x=\int_{u=1}^{u=z} \frac{z}{u^2}\left( \LL - \left(\frac{z}{u}\right)^{\LL^m-1}\right) \,du =\frac{\LL^{m+1}(z-1)-z^{\LL^m}+1}{\LL^m}.$$
We thus find that $\phi:=\Phi|_{y=1}=z-1$ satisfies  
\beq\label{f15}
x=\frac{\LL^{m+1}\phi-(1+\phi)^{\LL^m}+1}{\LL^m(\LL-1)}, \quad \psi:=\Psi|_{y=1}=(1+\phi)^{\HH_Y}.
\eeq 
When $\HH$ is the Hodge-Deligne polynomial, we have 
\beq\label{f20} \psi(x)=1+\sum_{n\ge 1}\frac{x^n}{n!}\HH_{Y[n]}.\eeq 
Therefore, we proved the following.
\begin{theorem} \cite[(4.22) and (4.24)]{Manin1999}
    The Hodge polynomials and the Poincar\'e polynomials $\HH_{Y[n]}$ of $Y[n]$ are uniquely determined by \eqref{f15} and \eqref{f20}.
\end{theorem}

\bigskip

\bibliographystyle{amsalpha}
\bibliography{references}

\end{document}